\documentclass[11pt]{amsart}
\usepackage{amsmath}
\usepackage{amsfonts}
\usepackage{amssymb}
\usepackage{latexsym}
\usepackage{amsthm}

\newtheorem{theorem}{Theorem}

\newcounter{other}            
\newtheorem{otherth}[other]{Theorem}              
\newtheorem{otherl}[other]{Lemma}        

\newcommand{\Cn}{\mathbb{C}^n}
\newcommand{\Sn}{\mathbb{S}_n}
\newcommand{\Bn}{\mathbb{B}_n}

\newcommand{\D}{\mathbb D}
\newcommand{\C}{\mathbb C}
\newcommand{\T}{\mathbb T}
\newcommand{\N}{\mathbb N}

\newcommand{\eps}{\varepsilon}

\newcommand{\intcirc}{\int_0^{2\pi}}

\newcommand{\Hi}{H^\infty}
\newcommand{\Li}{L^\infty}
\newcommand{\Omep}{\Omega_\varepsilon}
\newcommand{\Bl}{\mathcal B}

\begin{document}

\title{Closure of Hardy spaces in the Bloch space}

\author[P. Galanopoulos]{Petros Galanopoulos}
\address{Petros Galanopoulos \\School of Mathematics\\
\\ Aristotle University of Thessaloniki \\
54124
Thessaloniki \\
Greece} \email{petrosgala@math.auth.gr}

\author[N. Monreal Gal\'{a}n]{Nacho Monreal Gal\'{a}n}
\address{Nacho Monreal Gal\'{a}n, Departament of Mathematics, University of Crete, Voutes Campus, 70013 Heraklion, Crete, Greece.}
\email{nacho.mgalan@gmail.com}

\author[J. Pau]{Jordi Pau}
\address{Jordi Pau \\Departament de Matem\`{a}tica Aplicada i Analisi\\
Universitat de Barcelona\\
08007 Barcelona\\
Catalonia, Spain} \email{jordi.pau@ub.edu}



\thanks{ N.\ Monreal was supported in part by the project MTM2011-24606, and in part also by the research project PE1(3378) implemented within the framework of the Action ``Supporting Postdoctoral Researchers'' of the Operational Program ``Education and Lifelong Learning'' (Action's Beneficiary: General Secretariat for Research and Technology), co-financed by the European Social Fund (ESF) and the Greek State. Jordi Pau was
 supported by the DGICYT grant MTM$2011$-$27932$-$C02$-$01$
(MCyT/MEC)}

\begin{abstract}
A description of the Bloch functions that can be approximated in the Bloch norm by functions in the Hardy space $H^p$ of the unit ball of $\Cn$ for $0<p<\infty$ is given. When $0<p\leq1$, the result is new even in the case of the unit disk.
\end{abstract}

\maketitle




\section{Introduction.}

\noindent Let $\D$ and $\T$ be, respectively, the unit disk and the unit circle of the complex plane $\C$. For $0<p<\infty$, recall that the Hardy space $H^p(\D)$ is the space of
analytic functions $f$ in the unit disc such that
$$\|f\|^p_{p}=\sup_{0<r<1}\intcirc |f(re^{i\theta})|^p\, \frac{d\theta}{2\pi}<+\infty.$$
For $p=\infty$, $\Hi(\D)$ is the space of all bounded
analytic functions in the unit disk. Recall also that the Bloch space $\Bl(\D)$ is formed by the analytic functions $f$ on $\D$ such that
$$\|f\|_{\Bl}=\sup_{z\in\D}\,(1-|z|^2)\,|f'(z)|<\infty.$$

\noindent In \cite{MN}, a characterization of the closure in the Bloch norm of $H^p\cap\Bl$ for $1<p<\infty$ was given in terms of the area of certain non-tangential level sets of the Bloch function: given a function $f\in\Bl$ and
$\varepsilon>0$ define the level set of $f$ as
$$\Omep (f):=\{z\in\D\ :\ (1-|z|^2)|f'(z)|\geq\varepsilon\}.$$
\noindent Recall that a Stolz angle with vertex in $\zeta\in\T$ is the set
$$\Gamma(\zeta)=\Gamma_\alpha(\zeta):=\{z\in\D:|z-\zeta|<\dfrac{\alpha}{2}(1-|z|)\},$$
with $\alpha>2$, and that
$$A_h(\Omega):=\int_\Omega\frac{dA(z)}{(1-|z|^2)^2},$$
where $dA(z)$ is the area measure in $\D$, represents the hyperbolic area of $\Omega\subset\D$. Then the result is the following:

\begin{otherth}\label{thMN}
Let $f$ be a function in the Bloch space $\Bl$ and $1<p<\infty$. Then $f$ is in
the closure in the Bloch norm of $\mathcal{B}\cap H^p$ if and only if for any $\eps>0$ the function
$A_ h(\Gamma (\zeta)\cap \Omep(f))^{1/2}$ is in  $L^p(\T)$.
\end{otherth}

\noindent A basic tool in the proof of this result was the characterization of Hardy spaces in terms of the area function (a result due to Marcinkiewicz and A. Zygmund \cite{MZ} for $p>1$, and extended to the case $0<p\le 1$ by A. Calder\'{o}n \cite{Cal}), that is, for $0<p<\infty$, a function $f$ is in $H^p$ if and only if its corresponding Lusin Area function
$$A(f)(\zeta)=\left(\int_{\Gamma(\zeta)}|f'(z)|^2dA(z)\right)^{1/2}$$
is in $L^p(\T)$. The proof of Theorem \ref{thMN} was based on a previous result by P.\ Jones on the closure of $BMOA$ in $\Bl$ (see \cite{GZ}). The duality argument given in the proof in \cite{MN} can not be used for $0<p\leq 1$, so that this case requires of new techniques. In this paper we solve the case $0<p\le 1$. It turns out that the proof given works equally for all $0<p<\infty$, and furthermore, it may be done in the open unit ball $\Bn$ of the $n$-dimensional complex space $\Cn$. The case $p=\infty$ is still an open problem, and will be discussed in the last Section.

\bigskip

\noindent Now we are going to introduce some notation. For $z,w\in \Cn$, let
\[ \langle z,w \rangle =z_ 1 \bar{w}_ 1+\dots +z_ n \bar{w}_ n.\]
Hence, $|z|^2=\langle z,z\rangle$. In this context, for $0<p<\infty$ the Hardy space $H^p(\Bn)$ consists of those holomorphic functions $f$ on $\Bn$ such that
$$\|f\|^p_{p}=\sup_{0<r<1}\int_{\Sn} |f(r\zeta)|^p \, d\sigma(\zeta)<+\infty,$$
where $\Sn$ denotes the unit sphere in $\Cn$ and $\sigma$ is the normalized surface measure on $\Sn$. As in the case for $n=1$, for $p=\infty$ the corresponding space $\Hi(\Bn)$ is the space of bounded holomorphic functions defined on $\Bn$.

\noindent The Hardy space $H^p(\Bn)$ may be also characterized by means of a corresponding area function. In order to define it, let $Rf$ denote the radial derivative of $f$, that is,
\[Rf(z)= \sum_{k=1}^{n} z_ k \frac{\partial f}{\partial z_ k} (z),\qquad z=(z_ 1,\dots,z_ n)\in \Bn.\]
\noindent Besides, the hyperbolic measure in $\Bn$ is given by
$$d\lambda_ n(z)=\frac{dv(z)}{(1-|z|^2)^{n+1}},$$
where $dv(z)$ is the normalized volume measure in $\Cn$. The admissible Area function is then defined as
\begin{displaymath}
Af(\zeta)=\left ( \int_{\Gamma (\zeta)} |Rf(z)|^2 \,(1-|z|^2)^2\, d\lambda_ n(z)\right )^{1/2},\qquad \zeta \in \Sn,
\end{displaymath}
where $\Gamma(\zeta)$ denotes now the admissible Koranyi region, that is,
\begin{displaymath}
\Gamma(\zeta)=\Gamma_{\alpha}(\zeta):=\left \{z\in \Bn: |1-\langle z,\zeta\rangle |<\frac{\alpha}{2} (1-|z|^2) \right \}.
\end{displaymath}
When $n=1$ this region coincides with the usual Stolz angle in $\D$. The following result is the generalization of the area theorem for $\Bn$ and can be found, for example, in \cite{FS} or \cite[Theorem 5.3]{P1}.

\begin{otherth}\label{CMZ}
Let $0<p<\infty$, then $f\in H^p(\Bn)$ if and only if $Af\in L^p(\Sn)$. Furthermore, if $f(0)=0$ the norms $\|f\|_{H^p}$ and $\|Af\|_{L^p}$ are comparable.
\end{otherth}

\medskip

\noindent The Bloch space $\Bl:=\Bl(\Bn)$ consists of those functions $f$ holomorphic on $\Bn$ such that
\begin{displaymath}
\|f\|_{\mathcal{B}}=\sup_{z\in \Bn} (1-|z|^2) |Rf(z)| <\infty,
\end{displaymath}
\noindent This seminorm is not conformally invariant (for more details, see \cite[Chapter 3]{ZhuBn}), but it is equivalent to the seminorm defined above for the unidimensional case, and more convenient for the statements here. As in the one-dimensional case, $\Hi(\Bn)\subset\Bl(\Bn)$.

\smallskip

\noindent The level sets here are defined as
$$\Omep (f):=\{z\in\Bn\ :\ (1-|z|^2)|Rf(z)|\geq\varepsilon\}.$$
\noindent Clearly, the hyperbolic volume of any set $\Omega\subset\Bn$ is
\begin{displaymath}
V_ h(\Omega)=\int_{\Omega} d\lambda_ n.
\end{displaymath}

\medskip

\noindent The result proved in this paper is the following one.

\begin{theorem}\label{mainth}
Let $f$ be a function in the Bloch space $\Bl(\Bn)$ and $0<p<\infty$. Then $f$ is in
the closure in the Bloch norm of $\Bl\cap H^p(\Bn)$ if and only if for any $\eps>0$ the function
$V_ h(\Gamma (\zeta)\cap \Omep(f))^{1/2}$ is in  $L^p(\Sn)$.
\end{theorem}

\medskip

\noindent The necessity is done in the same way as in Theorem \ref{thMN}. The sufficiency is slightly harder. From \cite[p.51]{ZhuBn} one may express the function $f\in\Bl$ as
\begin{equation}\label{decomp}
f(z)=f(0)+\int_{\Bn} Rf(w)\,L(z,w) \,dv_{\beta}(w),
\end{equation}
where the kernel
\begin{displaymath}
 L(z,w)=\int_{0}^{1} \left ( \frac{1}{(1-t\langle z,w \rangle )^{n+1+\beta}}-1\right ) \,\frac{dt}{t}
\end{displaymath}
satisfies
\begin{displaymath}
 |L(z,w)| \le \frac{C}{|1-\langle z,w \rangle |^{n+\beta}} .
\end{displaymath}
\noindent Here
$$ dv_{\beta}(z)=c_{\beta}(1-|z|^2)^{\beta} \,dv(z),$$
where $\beta>-1$ and $c_{\beta}$ is a normalizing constant taken so that $v_{\beta}(\Bn)=1$.

\medskip

\noindent The following integral estimate has become indispensable in this area of Analysis. One may find the proof in \cite[Theorem 1.12]{ZhuBn}.

\begin{otherl}\label{IctBn}
Let $t>-1$ and $s>0$. There is a positive constant $C$ such that
\[ \int_{\Bn} \frac{(1-|w|^2)^t\,dv(w)}{|1-\langle z,w\rangle |^{n+1+t+s}}\le C\,(1-|z|^2)^{-s}\]
for all $z\in \Bn$.
\end{otherl}

\medskip

\section{Proof of Theorem \ref{mainth}.}

\begin{proof}
Let first $f$ be in the closure in the Bloch norm of $H^p\cap\Bl$. Then, given $\eps>0$ there exists $g\in H^p$ such that $\|f-g\|_{\Bl}<\eps/2$. As in \cite{MN}, one just needs to observe that $\Omep(f)\subseteq\Omega_{\eps/2}(g)$ to see that for any $\zeta\in\Sn$
$$V_h(\Gamma(\zeta)\cap\Omep(f))\leq\frac{4}{\eps}Ag(\zeta)^2.$$
\smallskip

\noindent Since $Ag\in L^p(\Sn)$, the necessity is proved.

\medskip

\noindent From now on, set $\Omep=\Omep(f)$. Let $f$ be in the Bloch space and assume that $V_ h(\Gamma (\zeta)\cap \Omep)^{1/2}$ is in  $L^p(\Sn)$. Given $\eps>0$, one wants to find $f_2\in\Bl\cap H^p$ such that $\|f-f_2\|_{\Bl}\leq\eps$. To this end, and applying (\ref{decomp}), $f$ may be decomposed in the sum of two functions $f(z)=f_ 1(z)+f_ 2(z)$, where
\begin{displaymath}
f_ 1(z)=\int_{\Bn\setminus \Omep} Rf(w)\,L(z,w) \,dv_{\beta}(w)
\end{displaymath}
\noindent and
\begin{displaymath}
f_ 2(z)=f(0)+\int_{\Omep} Rf(w)\,L(z,w) \,dv_{\beta}(w),
\end{displaymath}
with $\beta$ big enough to be fixed later. It is easy to see that
\begin{displaymath}
 \begin{split}
 |Rf_ 1(z)| &\leq \int_{\Bn\setminus \Omep} |Rf(w)|\left|\int_0^1\frac{\langle z,w\rangle}{1-t\langle z,w\rangle}dt\right|\,dv_{\beta}(w)\,\leq \\ &\leq \eps\,C(n,\beta)\,\int_{\Bn}\frac{(1-|w|^2)^{\beta-1}}{|1-\langle z,w\rangle|^{n+\beta+1}}dv(w).
 \end{split}
\end{displaymath}
\noindent Now Lemma \ref{IctBn} with $t=\beta-1$ and $s=1$ shows that $\|f_1\|_{\Bl}\lesssim\eps$.

\medskip

\noindent Thus it remains to see that $f_ 2 \in H^p(\Bn)$, or by Theorem \ref{CMZ}, that $Af_2\in L^p(\Sn)$. For this aim, one needs the following two Lemmas. The first one may be thought as a generalized version of Lemma \ref{IctBn}, and appears in \cite[Lemma 2.5]{of}.

\begin{otherl}\label{LI2-Bn}
Let $s>-1$, $r,t>0$, and $r+t-s>n+1$. If $t,r<s+n+1$ then, for $a,z\in \Bn$, one has
\begin{displaymath}
\int_{\Bn}\frac{(1-|w|^2)^s}{|1-\langle z,w \rangle|^r\,|1-\langle a,w \rangle |^t}\,dv(w)\leq
C\,\frac{1}{|1-\langle z,a \rangle |^{r+t-s-n-1}}.
\end{displaymath}
\end{otherl}

\medskip

\noindent The following estimation  may be found in \cite{Ar} and \cite{Jev}, and it is the analogue in $\Bn$ of \cite[Proposition 1]{Lue1}.

\begin{otherl}\label{AJL}
 Let $0<s<\infty$ and $b>n\max(1,1/s)$. Then there exists a constant $C>0$ depending on $s$, $b$ and on the angle of the region $\Gamma(\zeta)$ such that
 \begin{displaymath}
  \int_{\Sn} \left ( \int_{\Bn} \Big (\frac{1-|z|^2}{|1-\langle z,\zeta\rangle|}\Big )^b d\mu(z)\right )^s \,d\sigma(\zeta)\le C \int_{\Sn} \mu(\Gamma(\zeta))^s\,d\sigma(\zeta),
 \end{displaymath}
where $\mu$ is a positive measure on $\Bn$.
\end{otherl}

\medskip

\noindent Back to the proof, as in the previous estimation for $Rf_1$ one has that

\begin{displaymath}
\begin{split}
|Rf_ 2(z)|^2 &\le C(n,\beta) \,\|f\|^2_{\Bl} \left (\int_{\Omep} \frac{(1-|w|^2)^{\beta-1}}{|1-\langle z,w\rangle|^{n+\beta+1}}\,dv(w)\right )^{2}
\\
&\le \tilde{C}(n,\beta) \,\|f\|^2_{\Bl} \,\left (\int_{\Omep} \frac{(1-|w|^2)^{\beta-1}}{|1-\langle z,w\rangle|^{n+\beta+1}}\,dv(w)\right ) \,(1-|z|^2)^{-1},
\end{split}
\end{displaymath}
after an application of Lemma \ref{IctBn}. Then Fubini's theorem gives
\begin{displaymath}
\begin{split}
Af_ 2(\zeta)^2 &\lesssim \int_{\Gamma(\zeta)} \left (\int_{\Omep} \frac{(1-|w|^2)^{\beta-1}}{|1-\langle z,w\rangle|^{n+\beta+1}}\,dv(w)\right )\,\frac{dv(z)}{(1-|z|^2)^n}
\\
&=\int_{\Omep} \left (\int_{\Gamma(\zeta)} \frac{dv(z)}{(1-|z|^2)^n\,|1-\langle z,w\rangle|^{n+\beta+1}}\right )(1-|w|^2)^{\beta-1}dv(w).
\end{split}
\end{displaymath}
\noindent Now, $z\in\Gamma(\zeta)$ implies that $(1-|z|^2)\simeq |1-\langle z,\zeta\rangle|$. Then, applying Lemma \ref{LI2-Bn} with $r=t=n+\beta+1$ and $s=\beta+1$ one has that
\begin{displaymath}
\begin{split}
\int_{\Gamma(\zeta)} \frac{dv(z)}{(1-|z|^2)^n\,|1-\langle z,w\rangle|^{n+\beta+1}}
&\simeq \int_{\Gamma(\zeta)} \frac{(1-|z|^2)^{\beta+1}\,dv(z)}{|1-\langle z,\zeta\rangle|^{n+\beta+1}\,|1-\langle z,w\rangle|^{n+\beta+1}}
\\
 &\leq \int_{\Bn} \frac{(1-|z|^2)^{\beta+1}\,dv(z)}{|1-\langle z,\zeta\rangle|^{n+\beta+1}\,|1-\langle z,w\rangle|^{n+\beta+1}}
 \\
 &\lesssim \frac{1}{|1-\langle w,\zeta\rangle|^{n+\beta}}.
\end{split}
\end{displaymath}
Hence,
\begin{displaymath}
\begin{split}
\|Af_ 2 \|^p_{L^p} &\lesssim \int_{\Sn} \left ( \int_{\Omep} \frac{1}{|1-\langle w,\zeta\rangle|^{n+\beta}}(1-|w|^2)^{\beta-1}\,dv(w)\right )^{p/2} d\sigma(\zeta)
\\
&=\int_{\Sn} \left ( \int_{\Omep} \frac{(1-|w|^2)^{n+\beta}}{|1-\langle w,\zeta\rangle|^{n+\beta}}\frac{dv(w)}{(1-|w|^2)^{n+1}}\right)^{p/2} d\sigma(\zeta)
\\
&=\int_{\Sn} \left ( \int_{\Bn} \left(\frac{(1-|w|^2)}{|1-\langle w,\zeta\rangle|}\right)^{n+\beta}\,d\mu(w)\right)^{p/2} d\sigma(\zeta).
\end{split}
\end{displaymath}
Here
\[ d\mu (w)=\frac{\chi_{\Omep}(w)dv(w)}{(1-|w|^2)^{n+1}},\]
where $\chi_{\Omep}$ denotes the characteristic function of $\Omep$, is a positive Borel measure. Then Lemma \ref{AJL} with $s=p/2$ and $b=n+\beta$, where $\beta$ is positive and bigger than $n\cdot(2/p-1)$, shows that $Af_2\in L^p(\Sn)$. This finishes the proof.

\end{proof}

\section{The case $p=\infty$.}

\noindent From now on let $n=1$. The problem of describing the closure of the space of bounded analytic functions in the Bloch norm was posed in \cite{ACP}, and still remains open. Remember that $\Hi\subset\Bl$. Theorem \ref{CMZ} does not hold for $p=\infty$, so the proof given here does not work in this case. Nevertheless, it is interesting to outline that the proof given holds also if one considers the class of analytic functions with area function in $\Li$.

\medskip

\noindent One may see, instead, that the analogue for $p=\infty$ of the condition given in Theorem \ref{mainth} (that is, $A_h(\Omep(f)\cap\Gamma(\zeta))\in \Li(\T)$) is not necessary for a function $f$ to be in the closure in the Bloch norm of the space of bounded analytic functions. To this end, for $k\in\N$ take the points $z_k=1-2^{-k}$, and consider the sequence $\{z_k\}$, which is a radial separated sequence. In particular, $\{z_ k\}$ is an interpolating sequence for $H^{\infty}$ (see \cite[Chapter VII, p.279]{Gar}). By Carleson interpolation theorem, there exists $\delta>0$ such that
$$(1-|z_k|^2)|B'(z_k)|\geq\delta,$$
where  $B$ denotes the Blaschke product with zeros $\{z_k\}$, which is clearly in $\Hi$.
\noindent Given now $\eps<\delta/4$ there exists $\rho>0$ such that $(1-|z|^2)|B'(z)|\geq\eps$ on each $D_h(z_k,\rho)$, that is, the hyperbolic disk with center $z_k$ and radius $\rho$. Hence,
$$\bigcup_{k\in\N}D_h(z_k,\rho)\subset\Omep(B)\cap \Gamma (1).$$
\noindent Since the sequence is separated, one can take $\rho>0$ so that the disks $D_h(z_k,\rho)$ are pairwise disjoints. Now it is easy to see that $A_h(\Omep(B)\cap\Gamma(\zeta))$ is not in $\Li(\T)$, since $A_h(D_h(z_k,\rho))\geq C$ for a certain constant $C>0$ only depending on $\rho$.

\medskip

\noindent In \cite[Section 3.6]{Xiao} one may find a sufficient condition for a Bloch function to be in the closure in the Bloch norm of $\Hi(\D)$. The condition is the following: For every $\eps>0$ one has that
\begin{equation}\label{xiao_conj}
 \sup_{w\in\D}\int_{\Omep(f)}\frac{1}{|1-\overline{w}z|^2}\,dA(z)<\infty.
\end{equation}

\noindent The sufficiency of this condition is checked following also the proof in \cite{MN}. Let $f\in\Bl$ satisfying (\ref{xiao_conj}). Without loss of generality one may take $f(0)=f'(0)=0$. Hence the function can be expressed by the following integral (see \cite[Proposition 4.27]{Zhu})
\begin{displaymath}
 f(w)=\int_\D\frac{(1-|z|^2)f'(z)}{\overline{z}(1-\overline{z}w)^2}\,dA(z)=f_1(w)+f_2(w),
\end{displaymath}
\noindent where $f_1$ and $f_2$ are taken as in the previous Section. Then one may see that $\|f_1\|_{\Bl}\lesssim\eps$. In order to see that $f_2\in\Hi$ one just have to observe that
\begin{displaymath}
 |f_2(w)|\leq\int_{\Omep(f)}\frac{(1-|z|^2)|f'(z)|}{|\overline{z}||1-\overline{w}z|^2}\,dA(z)\leq C(\eps)\, \|f\|_{\Bl}\int_{\Omep(f)}\frac{dA(z)}{|1-\overline{w}z|^2},
\end{displaymath}
\noindent and then apply the hypothesis.

\medskip

\noindent Actually, in \cite[p.71]{Xiao}, J. Xiao  conjectured that condition (\ref{xiao_conj}) is also necessary. Nevertheless, the same example as above gives a counterexample to that conjecture, just by evaluating the integral for $w$ approaching 1 non-tangentially. Indeed, if $w_ m=1-2^{-m}$, then
\[
\int_{\Omep(B)}\frac{1}{|1-\overline{w}_ m z|^2}\,dA(z) \ge \sum_{k=1}^{m}\int_{D_ h(z_ k,\rho)}\frac{1}{|1-\overline{w}_ m z|^2}\,dA(z).
\]
Now, it is easy to see that there is a constant $C$ depending on $\rho$ such that $|1-\overline{w}_ m z|\le C (1-|z_ k|)$ for $z\in D_ h(z_ k,\rho)$ and $k\le m$. This clearly implies that 
\[
\int_{\Omep(B)}\frac{1}{|1-\overline{w}_ m z|^2}\,dA(z) \longrightarrow +\infty
\]
as $m\rightarrow \infty$, proving that the condition (\ref{xiao_conj}) is not necessary.


\end{document}